\documentclass{article}
\usepackage{amsmath}
\usepackage{amsfonts}

\def\a{\alpha}
\def\b{\beta}

\def\d{\downarrow}
\def\N{\mathbb{N}}
\def\R{\mathbb{R}}
\def\C{{\cal C}}

\def\n{\noindent}

\def\p{\partial}
\def\r{\rightarrow}

\def\O#1{{\bar{#1}}}

\newtheorem{cor}{Corollary}
\newtheorem{lem}{Lemma}
\newtheorem{prop}{Proposition}
\newtheorem{theo}{Theorem}
\newtheorem{obs}{Observation}

\newcommand{\be}{\begin{equation}}
\newcommand{\ee}{\end{equation}}

\title{Small growth vectors of the compactifications\\
of the contact systems on $J^{\,r}(1,1)$}

\author{Piotr Mormul\\
          \\
\small Mathematical Institute, Polish Academy of Sciences,\\
\small {\'S}niadeckich \,str. 8, \,00-097 Warsaw, \,Poland.\\
\small \\
\small On leave from Institute of Mathematics, University of Warsaw,\\
\small e-mail: mormul@mimuw.edu.pl\normalsize}
\date{}

\begin{document}
\maketitle

%%%%%%%%%%%%%%%%%%%%%%%
\begin{abstract}
\n It is well known that the compactifications of the canonical contact 
systems living on real jet spaces $J^{\,r}(1,1)$, $r \ge 2$, are locally 
universal Goursat distributions, $\Delta^r$, living on compact manifolds 
(called Goursat monsters) having open dense jet-like ($J^{\,r}(1,1)$-like) 
parts. 

\n By virtue of the results of Jean (1996), one was able, for each 
$r \ge 2$, to recursively compute the {\it small growth vector\,} of 
\,$\Delta^r$ at any point of the $r$-th monster. The result was got by 
performing series of $r$ operations taken, in function of the local 
geometry of \,$\Delta^r$ in question, from the set of fixed recursive rules 
$\{\,1,\,2,\,3\,\}$ (called in the present text S,\,T,\,G, respectively). 
By the local universality of \,$\Delta^r$ one was thus able to compute 
{\it all\,} small growth vectors of all existing Goursat distributions. 

\n In the work of Mormul (2004) proposed were explicit solutions to 
the series of Jean's recurrences. The solutions uncovered a surprisingly 
involved underlying arithmetics -- a kind of G{\"o}del-like encoding of 
words over a three-letter alphabet $\{$G,\,S,\,T$\}$ by neat sequences 
of positive integers. Yet, those formulas, though characterizing really 
existing objects, appeared as if from thin air, and proofs were postponed 
to another publication. 

\n In the present contribution we submit proofs of our formulas from 2004. 
It is not, however, a plain check that our candidates satisfy Jean's 
recurrences. Under the way we retrieve and re-produce those surprising 
formulas. 
\vskip2.5mm
\n{\it Key words and phrases.} Goursat flag, small growth vector, 
solution of Jean's recurrences, encoding of letter words by 
sequences of integers
\vskip2.5mm
\n 2000 {\it Mathematics Subject Classification.} 58A17,\,\,30.
\end{abstract}
%%%%%%%%%%%%%
%%%%%%%%%%%%%%
%%%%%%%%%%%%%%%
\section{Goursat distributions and their small growth vectors}
\label{vectors}
Goursat flags are certain special nested sequences, in general of 
variable length $r$ ($2 \le r \le n - 2$) of subbundles in the tangent 
bundle $TM$ to a smooth or analytic $n$-dimensional manifold $M$: 
$D^r \subset D^{r-1} \subset \cdots \subset D^1 \subset D^0 = TM$. 
Namely, one demands, for $l = r,\,r-1,\dots,1$, that (a)\,the Lie 
square of $D^l$ be $D^{l-1}$ and (b)\,$\text{cork}\,D^l = l$. In 
other words, that the {\it big growth vector\,} of the distribution 
$D^r$ be $[n-r,\,n-r+1,\,n-r+2,\dots,\,n-1,\,n]$ at each point of $M$. 
\vskip1.2mm
This, very restricted, class of objects was being investigated 
(intermittently) over the last 120 years, with contributions, among 
others, by E.\,von\,Weber \cite{vW} and E.\,Cartan \cite{C}. They 
independently proved that every corank-$r$ Goursat distribution $D^r$ 
on $M$ around a generic point locally behaves in a unique way visualised 
by the {\it chained model} -- the germ at $0 \in \R^n(x^1,\dots,x^{r+2};
\,x^{r+3},\dots,x^n)$ of 
%%%%%%
\be
\big(\,\p_n,\dots,\,\p_{r+3};\;\p_{r+2},\;\,\p_1 + x^3\p_2 + x^4\p_3 + \cdots + x^{r+2}\p_{r+1}\,\big)\tag{A}
\ee
(these are vector fields generators; effectively used are only first $r + 2$ 
coordinates). Chained models can be viewed as the simplest instance of a family 
of local writings (preliminary normal forms) of Goursat distributions, obtained 
much later by Kumpera and Ruiz. As a matter of fact, Kumpera and Ruiz discovered, 
to a big surprise of the mathematical community, {\it singularities\,} hidden in 
flags, and their pseudo-normal forms were merely a by-product -- visualisations 
of those singularities. 

When $n = r + 2$, a Goursat distribution $D^r$ on $M$ has around generic 
points a unique local description as the germ at $0 \in \R^{r+2}(t,\,x,\,x_1,\dots,\,x_r)$ 
of 
%%%%%%
\be
\big(\,\p/\p{x_r},\;\,\p/\p t + x_1\p/\p x + x_2\p/\p x_1 + \cdots + x_r\p/\p x_{r-1}\,\big)\,.\tag{B}
\ee
That is, of the canonical contact system on $J^{\,r}(1,1)$. (The independent 
variable is $t$, the dependent one is $x$, and $x_1,\,x_2,\dots,\,x_r$ 
can be viewed as the consecutive derivatives of $x$ with respect to $t$.) 
This is precisely the result obtained independently by von\,Weber and Cartan. 
\vskip1.2mm
%%%%%%%%%%%%%%%%%%%%%
For a theory of Goursat flags it has been essential to have a flexible 
length $r$ that could be increased without changing a manifold, allowing 
for flags to be prolonged in length on one and the same manifold. In what 
follows, however, we deal only with the most classical situation $n = r + 2$, 
when the initial Goursat distribution $D^r$ is of rank 2 and its flag is 
already of maximal length $r = n - 2$. 
For, as long as Goursat distributions are considered locally, this causes 
no loss of generality. In fact, when $n > r + 2$, in a Goursat $D^r$ 
locally there always splits off an integrable subdistribution of rank 
$n - r - 2$ which leaves a second direct summand (of rank 2) inside $D^r$ 
invariant. Like $\big(\p_n,\dots,\p_{r+3}\big)$ in (A) above, which leaves 
invariant the summand 
$\big(\,\p_{r+2},\;\p_1 + x^3\p_2 + x^4\p_3 + \cdots + x^{r+2}\p_{r+1}\,\big)$. 
Compare, for instance, Corollary 1.3 in \cite{arithm}. 
\vskip1.2mm
Note that the rank-2 and corank-$r$ Goursat distributions have the big 
growth vector $[2,\,3,\,4,\dots,\,r+1,\,r+2]$ at each point. 
\vskip1.2mm
%%%%%%%%D-1
\n{\bf Definition 1.} 
Small growth vector (sgrv in all what follows) $\{n_j(p)\}$ of a distribution 
$D$ at a point $p \in M$ is the sequence of dimensions at $p$ of the (local) 
modules of vector fields $V_j$, \,$V_1 = D$, \,$V_{j+1} = V_j + [D,\,V_j]$,
$$
n_j(p) = {\rm dim}\,V_j(p)\,,\qquad j = 1,\,2,\,3,\dots
$$
We are interested only in distributions for which this sequence, for 
every $p \in M$, attains (sooner or later) the value ${\rm dim}\,M$, 
with `sooner or later' underscored. Such distributions are called 
{\it completely nonholonomic}. When, for such a distribution $D$, 
there happens $n_{l-1}(p) < n_l(p) = {\rm dim}\,M$, then $l$ is 
called the {\it nonholonomy degree\,} of $D$ at $p$. 
\vskip1.5mm
%%%%%%%%%%%%%
It is a short Lie algebra exercise to show that the Goursat distributions 
are completely nonholonomic. Around jet-like points $p$ where the visualisation 
(B) is effective, there holds $n_j(p) = 1 + j$ for $j = 1, \,2,\dots,\,n - 1$, 
and then the nonholonomy degree, minimal possible for Goursat in dimension $n$, 
is $n - 1 = r + 1$, the number $r$ being the order of jets in the visualisation. 
For non-jet-like points, even {\it much\,} slower sgrv's (and hence bigger 
nonholonomy degrees) are possible, as is invoiced in \cite{arithm} and will 
be eventually clear from the present work. The biggest nonholonomy degree 
in dimension $n$ will turn out to be $F_n$, the $n$-th Fibonacci number. 
It is about $\frac{1}{\sqrt{5}}\left(\frac{1 + \sqrt{5}}{2}\right)^n$, 
compare Remark 2 below for more on that. 
\vskip1.5mm
%%%%%%%%%%%%%%%%%%%%%%%%%%%%%%%%%%%%%%%%
A given Goursat distribution $D^r$ may not feature all possible singularities 
in corank $r$, nor all possible sgrv's of Goursat distributions of that corank. 
But in each corank $r$ there do exist manifolds with Goursat structures of 
corank $r$ on them that are locally universal -- feature all possible singularities, 
hence also all possible sgrv's of Goursat distributions of corank $r$. 
Such are, for instance, the kinematical models extensively used and analysed 
in \cite{J}, and especially {\it monster Goursat manifolds} ${\mathbb P}^r(\R^2)$, 
with locally universal Goursat structures $\Delta^r$ living on them, constructed 
in \cite{MZ2}.\footnote{\,\,Goursat monsters could be alternatively built over, 
say, $S^2$ or $T^2$ or yet another compact 2-dimensional manifold, instead of 
$\R^2$. But they would serve just the same objective of local universality, 
while it is much easier to work with ${\mathbb P}^r(\R^2)$, notwithstanding 
$\R^2$ is not compact.} 
\vskip1.5mm
A natural question reads what are all occurring sgrv's of Goursat distributions 
of a fixed corank $r \ge 2$. Or, which is the same for that fixed $r$, what are 
all sgrv's of the distribution $\Delta^r$ living on ${\mathbb P}^r(\R^2)$. 
\vskip2mm
%%%%%%%%%%%%%%%%%%
We recapitulate here the information necessary for the present contribution. 
First of all, the germs of Goursat flags of length $r$ are partitioned into 
$F_{2r-3}$ invariant {\it geometric classes\,} labelled (encoded) by words 
of length $r$ over the alphabet $\{$G,\,S,\,T$\}$ starting with GG and such 
that GT is not allowed.\footnote{\,\,concerning the count of such words, 
yielding $F_{2r - 3}$, see Observation 1.10 in \cite{arithm}} We call such 
words {\it admissible}. The roots of this concept can be traced back to the 
paper \cite{J}, where the prototypes of geometric classes were called {\it regions}. 
By virtue of Jean's results \cite{J}, the sgrv of a Goursat germ depends solely 
on its geometric class. Hence, in length $r$ there are not more than $F_{2r-3}$ 
different sgrv's of Goursat germs. 
\vskip1.5mm
%%%%%%%%%%%%%%%%%%%%%%%%%%R-1
\n{\bf Remark 1.} It was not explained in \cite{J} whether the small growth 
vectors attached to the geometric classes, in any fixed length, were all different. 
It had remained a fine point for quite a time, having become completely clear 
only on the solutions' side of Jean's recurrences. See in this respect in 
\cite{arithm}: the discussion in Section 2.1 and Theorem 3.5.
\vskip1.1mm
Those who prefer the universal objects may think of the geometric classes 
in length $r$ as just the strata of a [very regular, besides] stratification 
of the $r$-th monster ${\mathbb P}^r(\R^2)$. Those strata are embedded 
submanifolds of ${\mathbb P}^r(\R^2)$ of codimensions equal to the numbers 
of non-G letters in their codes. (As a matter of fact, Jean's stratification 
of ${\mathbb P}^r(\R)^2$ is labelled in \cite{MZ2} by $r$-letter words over 
a newer alphabet $\{$R,\,V,\,T$\}$, with R replacing the previous G, V replacing 
S, and T replacing T.) \,Repeating, then, within each stratum, or geometric 
class, the sgrv of $\Delta^r$ is constant, not depending on points. 
\vskip1.5mm
In fact, it is a nondecreasing sequence of integers that starts with 2 and 
takes on, for the Goursat distributions, all intermediate integer values up 
to $r + 2$ inclusively. The key issue is to ascertain with what multiplicities 
do these integers appear. We recall that the sequence of multiplicities of 
integers in the sgrv of a given geometric class has been called in \cite{arithm} 
the {\bf derived vector} of that class. We stick to this terminology in 
the present paper. Sometimes, for brevity's sake, we will speak about the 
derived vector of an admissible word (= the code of a geometric class). 
\vskip2mm
%%%%%%%%%%%%%%%%%%%%%%%%%%%%%
The recipes to compute the derived vectors out of admissible (G,S,T) words 
which encode the geometric classes, were given in \cite{arithm}: firstly on 
a recursive basis rephrasing the results of \cite{J}, then in closed form 
formulas, which explicitly solved the recurrences for derived vectors. 
The first way did not need a proof, because Jean's recurrences, for the 
functions {\it beta\,} he used, automatically implied the recurrences for 
derived vectors being nothing but the sequences of {\it increments\,} of 
the beta functions. 

The second way necessitated a proof, or proofs of three separate statements 
in \cite{arithm}: Proposition 3.2, Theorem 3.3, and Theorem 3.4. The purpose 
of the present text is to furnish such proofs, after a slight reformatting 
of the statements under consideration: a part of old Proposition 3.2 and 
Theorem 3.3 are now merged into Theorem \ref{T-1}, while the rest of 
Proposition 3.2 and Theorem 3.4 are merged into Theorem \ref{T-2}. 
\vskip1.5mm
In \cite{arithm}, for each geometric class \,$\C$ of length $r$, the derived 
vector of \,$\C$, denoted by ${\rm der}(\C)$, is a function $\{2,\,3,\dots,\,r+1\} 
\r \N$, in which ${\rm der}(\C)(j)$ is the multiplicity of $j$ in the sgrv of 
\,$\C$, for every integer $2 \le j \le r+1$. Eventually the sgrv ends with 
the ceiling value $r + 2$. 
\vskip1mm
So how does Jean recursively arrive (modulo the mentioned passing to 
the increments of his original functions beta) at the functions 
${\rm der}(\C)$?
%%%%%%%%%%%
\begin{theo}[\cite{J}]\label{com-d}
For each geometric class \,$\C$ in length $r$, the function ${\rm der}(\C)$ 
equals the last term \,$d^{\,r}$ in the sequence of functions \,$d^{\,1},
\,d^{\,2},\dots,\,d^{\,r}$ \,that are constructed, by means of operations 
{\rm \,G,\,S,\,T} defined below, as follows. 
\end{theo}
One defines, recursively for $j = 1,\,2,\dots,\,r$, functions 
\,$d^{\,j} \colon \{2,\dots,\,j+1\} \r \N$. One starts by declaring 
$d^{\,1} = (1)$ and $d^{\,2} = (1,\,1)$, and then continues
$$
d^{\,j+2} = \begin{cases} {\rm G}\,(d^{\,j+1})\,, & \text{when the 
$(j+2)$-nd letter in the code of \,$\C$ is G},\\
                        {\rm S}\,(d^{\,j},\,d^{\,j+1})\,, & \text{when 
the $(j+2)$-nd letter in the code of \,$\C$ is \,S},\\
                        {\rm T}\,(d^{\,j},\,d^{\,j+1})\,, & \text{when 
the $(j+2)$-nd letter in the code of \,$\C$ is T}.
\end{cases}
$$
Concerning the operations in use, the simplest among them is the operation G,  
$$
{\rm G}(\a) \,= \left(\begin{array}{ccccc}
2 & 3 & 4 & 5 & \dots\\
\d & \d & \d & \d & \dots\\
1 & \a(2) & \a(3) & \a(4) & \dots
\end{array}\right). 
$$
So it is an "insert a 1 on the left and shift an argument vector $\a$ by 
one slot to the right" operation. Clearly, G($\a$) is a sequence (function) 
by one entry longer than $\a$ (with the domain by one element bigger than 
the domain of $\a$). 

\n The remaining operations S and T are more involved. They are two-argument, 
not just one-argument like G. Their first argument, say $\a$, is a sequence 
by one entry shorter than the second argument, say $\b$, and their outputs 
are yet by one entry longer than $\b$, 
$$
{\rm S}(\a,\,\b) \,= \left(\begin{array}{cccccc}
2 & 3 & 4 & 5 & 6 & \dots\\
\d & \d & \d & \d & \d & \dots\\
1 & 1 & \a(2) + \b(3) & \a(3) + \b(4) & \a(4) + \b(5) & \dots
\end{array}\right), 
$$
$$
{\rm T}(\a,\,\b) \,= \left(\begin{array}{cccccc}
2 & 3 & 4 & 5 & 6 & \dots\\
\d & \d & \d & \d & \d & \dots\\
1 & 1 & 2\b(3) - \a(2) & 2\b(4) - \a(3) & 2\b(5) - \a(4) & \dots
\end{array}\right). 
$$
So, in perhaps more catching terms, S is a {\it Fibonacci-like rule}, 
with appropriate shifts of the input vectors $\a$ and $\b$, while T 
is an {\it arithmetic progression rule}, with likewise shifts of the 
inputs. In fact, when the two inputs $\a$ and $\b$ are identically 
indented on the left, then their S (T, resp.) output starts, one 
row below, by inserting two 1's on the left and then performing 
the Fibonacci rule (arithmetic progression rule, resp.) on $\a$ 
and $\b$ skew-wise from the NW to the SE. To give an instance of 
each of these two-argument operations in action, and skipping 
brackets for bigger transparence,
$$
\begin{array}{rccccccc}
\a & = & 1 & 1 & 2 & 3 & & \\
\b & = & 1 & 1 & 1 & 2 & 3 & \\
{\rm S}(\a,\b) & = & {\bf 1} & {\bf 1} & 2 & 2 & 4 & 6
\end{array}
$$
and 
$$
\begin{array}{rcccccc}
\a & = & 1 & 1 & 2 & & \\
\b & = & 1 & 1 & 2 & 3 & \\
{\rm T}(\a,\b) & = & {\bf 1} & {\bf 1} & 1 & 3 & 4 
\end{array}
$$
(the two initial 1's inserted in the outputs on the left are written in bold). 
\vskip1.5mm
%%%%%%%%%%%%%%%R-2
\n{\bf Remark 2.} When, for any fixed length $r$, does the sum of values of 
${\rm der}(\C)$ attain the maximal value, or: when the nonholonomy degree 
$1 + \sum_{i=2}^{r+1}{\rm der}(\C)(i) = \max$\,? It is an exercise in the 
G, S, and T operations to show that, for each $r$, the maximal nonholonomy 
degree is realized only by the [most singular] class \,$\C = {\rm G\,G}
\underbrace{{\rm S\,S}\dots{\rm S}}_{r-2}$. \,Then \,${\rm der}(\C)(i) 
= F_{i-1}$ \,for \,$i = 2,\,3,\dots,\,r+1$, \,and 
$$
1 + \sum_{i=2}^{r+1}{\rm der}(\C)(i) \,= \,1 + F_1 + F_2 + \cdots + F_r \,= \,F_{r+2}\,,
$$
as invoiced earlier in this section. (Thus the nonholonomy degree of the distribution 
$\Delta^r$ always varies between $r + 1$ and $F_{r + 2}$ and these bounds are 
sharp for each $r \ge 2$. For instance, for $\Delta^3$, it varies between 4 and 
5. In general, however, {\it not\,} all intermediate integer values are realized 
by the nonholonomy degrees of $\Delta^r$. The smallest length with this happening 
is 6: the value 20 located between 7 and $21 = F_{6+2}$ is not realized, cf. 
Theorem 23 in \cite{CMor}.) 
%%%%%%%%%%%%%%%%
%%%%%%%%%%%%%%%%
\subsection{Useful aggregates of geometric classes.}
In what follows, the length ($r$) of a geometric class will not be essential. 
Much more important will be the number, say $s + 1 \ge 1$, of letters S in 
class' code that can otherwise be arbitrarily long.\footnote{\,\,By Theorem 
\ref{com-d}, in any fixed length $r \ge 2$, the generic jet-like geometric 
class with no letters S, written shortly ${\rm G}_r$, and only this class, 
has the simplest derived vector $(1,\,1,\dots,\,1)$ ($r$ 1's). Within ${\rm G}_r$, 
and only within it (in that length), the small growth vector of a Goursat 
distribution is the same as its big growth vector. In the paper we only 
consider the {\it remaining\,} geometric classes featuring certain letters 
S in their codes. That is, we consider the singular geometric classes.} 
To word our results, we arbitrarily {\bf fix this discrete parameter} 
$s \ge 0$. 
\vskip1mm
We have to introduce some more, still necessary, parameters -- the numbers
of letters T and G in classes' words going past the letters S, those letters 
S being run {\it backwards\,} from code's right end to left. 

\n Namely, agree that the last S in a code is followed by $k_0 \ge 0$ 
letters T, and then by $l_0 \ge 0$ letters G. Agree also that the one before 
last S is followed by $k_1 \ge 0$ letters T and then by $l_1 \ge 0$ letters G, 
and so on backwards until the first appearing letter S being followed by 
$k_s \ge 0$ letters T, then by $l_s \ge 0$ letters G. 

\n On top of that, let a code start with $l_{s+1} \ge 2$ letters G. 
(Except for this last quantity, all introduced integer parameters may 
even vanish, as it happens in the `Fibonacci' classes ${\rm GGS}_{s+1}$ 
discussed in Remark 2 above.) 
\vskip1.5mm
%%%%%%%%%%%%%%%%%%%%%%%%%R-3
\n{\bf Remark 3.} The indexation of integer parameters $k_j$ and $l_j$ proposed 
here slightly differs, by a backward shift of indices by one, from the indexation 
in \cite{arithm}. The same applies to the, central for the contribution, family 
(depending on the parameters $k_j$) of sequences $(A_0,\,A_1,\,A_2,\,\dots)$ 
introduced below before Theorem \ref{T-1}. The formulas, especially in 
Theorem \ref{T-2} below, assume thus a more compact view than in 
\cite{arithm}. 
\vskip1.5mm
In terms of these parameters, what can be said about the derived vectors 
emerging from Theorem \ref{com-d}\,?
\vskip2mm
%%%%%%%%%%%%%%%%%%%%
We will first answer the question for the geometric classes having 
$l_1 = l_2 = \dots = l_s = 0$. That is, for classes in whose codes the letters G 
occur only in the beginning (in number $l_{s+1} \ge 2$) and possibly in the end 
(in number $l_0 \ge 0$). The union of such geometric classes has been called 
in \cite{arithm} the hyperclass $0_s$. In the present text we use the concept 
of a hyperclass only implicitly, not explicitly. 

To word the relevant theorem, let us introduce a {\it family\,} of integer 
sequences 
%%%%%%%
\begin{itemize}
\item 
$A_0 = 1,\quad A_1 = 2 + k_0,\quad A_j = A_{j-2} + (1 + k_{j-1})A_{j-1}$\quad 
for\quad $2 \le j \le s+1$\,.
\end{itemize}
(This family is parametrized by the non-negative integers $k_0,\,k_1,\dots,\,k_s$.) 
%%%%%%%%%%%%%%%%%%%%%%%%%%%%%%%%%%%%%%
%%%%%%%%%%%%%%%%%%%%%%%%%%%%%%%%%%%%%%
%%%%%%%%%%%%%%%%%%%%%%%%%%%%%%%%%%%%%%
\section{Main theorems}
%%%%%%%%%%%%%
\begin{theo}[\cite{arithm}]\label{T-1}
In the derived vectors of the geometric classes with $s + 1$ letters {\rm \,S} 
and such that \,$l_1 = l_2 = \dots = l_s = 0$, there appear the following $s + 2$ 
different values, in the growing order: $A_0,\,A_1,\dots,\,A_{s+1}$. All these 
derived vectors are non-decreasing and, concerning the multiplicities of the 
listed values: 
\end{theo}
%%%%%%%%%%
\begin{itemize}
\item $A_0$ appears $2 + k_0 + l_0$ times in row; 
\item for $1 \le j \le s$, the number $A_j$ appears $1 + k_j$ times in row; 
\item $A_{s+1}$ appears $l_{s+1} - 1$ times in row. 
\end{itemize}
\vskip1mm
\n Proof of this theorem is given in Section \ref{pfT-1}. 
\vskip2.5mm
Before describing the sgrv's of the remaining geometric classes, we need an infinite 
series of sequences of the type $A$. This series will be parametrized by yet another 
natural number $N$. ($N$ is subject to some restrictions in function of $s$, but 
$s$ is not {\`a} priori bounded from above.)
%%%%%%%%%%%%%
\begin{itemize}
\item $A^{+N}_0 = 1,\quad A^{+N}_1 = 2 + k_N$\ \,($1 \le N \le s$)\,, 
\item $A^{+N}_j = A^{+N}_{j-2} + (1 + k_{j-1+N})A^{+N}_{j-1}$\quad\ for
\quad\ $2 \le j$, \,\,$N + j \le s + 1$\,.
\end{itemize}
We also have to specifically put in relief those letters S in a class' code 
that are followed by {\bf non-zero numbers} of letters G (going past some 
letters T, if any). They will play a central role in the characterization 
of the derived vector of such a class. 
\vskip1mm
Let, for a given geometric class not served by Theorem \ref{T-1} (hence, 
in particular, having $s + 1 \ge 2$ letters S in its code), 
$$
1 \le n_1 < n_2 < \cdots < n_q \le s
$$ 
be {\it all\,} indices $i$ between 1 and $s$ such that $l_i > 0$. (The index $n_1$ 
exists for each such class. The index $n_2$ -- already not for each such class. 
For its existence, the parameter $s$ should be at least 2 and, among the discrete 
parameters $l_1,\,l_2,\dots,\,l_s$ at least two should be non-zero, etc.) 
With these notations, we are now ready to complete the information given 
in Theorem \ref{T-1} by the following 
%%%%%%%%%%%%%%%%%%%%%%%%
\begin{theo}[\cite{arithm}]\label{T-2}
In the derived vectors of the geometric classes with $s + 1$ letters 
{\rm \,S} and such that, among the parameters $\{l_1,\,l_2,\dots,\,l_s\}$, exactly 
$l_{n_1},\,l_{n_2},\dots,\,l_{n_q}$ are non-zero (positive), $1 \le n_1 < n_2 < \cdots < n_q \le s$, 
\,$1 \le q \le s$, there appear $s + 2$ different values, listed below in 
the growing order, in the following $q + 1$ separate rows: 
%%%%%%%%%%%
\begin{align*}
&A_0,\,A_1,\,\dots,\,A_{n_1 - 1}\,; \\
%%%
A_{n_1}\big( & A_0^{+n_1},\;A_1^{+n_1},\,\dots,\;A_{n_2 - n_1 - 1}^{+n_1}\big)\,; \\
%%%
A_{n_1}A_{n_2 - n_1}^{+n_1}\big( & A_0^{+n_2},\;A_1^{+n_2},\,\dots,\;A_{n_3 - n_2 - 1}^{+n_2}\big)\,; \\
%%%
\dots \qquad & \dots \qquad \dots \qquad \dots \qquad \dots\\
%%%
A_{n_1}\prod_{j=1}^{q-2}A_{n_{j+1} - n_j}^{+n_j}\Big( & A_0^{+n_{q-1}},\;A_1^{+n_{q-1}},\,\dots,
\;A_{n_q - n_{q-1} - 1}^{+n_{q-1}}\Big)\,; \\
%%%
A_{n_1}\prod_{j=1}^{q-1}A_{n_{j+1} - n_j}^{+n_j}\Big( & A_0^{+n_q},\;A_1^{+n_q},\,\dots,
\;A_{s - n_q + 1}^{+n_q}\Big)\,.
\end{align*}
These derived vectors are always non-decreasing and, among the above-listed values: 
\end{theo}
%%%%%%%
\begin{itemize}
\item 
The $0$-th of them ($A_0$) appears $2 + k_0 + l_0$ times in row; 
\item
for $1 \le j \le s$, the $j$-th value appears $1 + k_j + l_j$ times in row 
(that is, $1 + k_j$ times, when $j \notin \{n_1,\,n_2,\dots,\,n_q\}$);
\item
the last and biggest $(s+1)$-st value appears $l_{s+1} - 1$ times in row. 
\end{itemize}
%%%%%%%%%%%%%%%%%%%%
\vskip1mm
\n Observations. (i)\,The table of $s + 2$ different values in ${\rm der}(\C)$ 
in Theorem \ref{T-2} is now formatted slightly differently than in Theorem 3.4 
in \cite{arithm}. The last entries in the rows (except the last row) {\it there\,} 
are {\it now\,} shifted to the front of the following rows. This is done by purely 
technical reasons of making the proof (in Section \ref{pfT-2}) more transparent. 
\vskip1mm
\n(ii)\,Also, to simplify the presentation in Theorem \ref{T-2}, the numbers 
standing {\it before\,} the brackets in the all but first rows in the above list 
do multiply all number entries standing {\it inside\,} the relevant brackets. 
\vskip1mm
\n(iii)\,When $n_1 = 1$, or when $n_{i+1} - n_i = 1$, then there is only one entry 
in the relevant row in the list. When $n_q = s$, there are only two entries in 
the last row. 
\vskip2.5mm
\n Proof of Theorem \ref{T-2} is given in Section \ref{pfT-2}.
\vskip2.5mm
%%%%%%%%%%%%%%%%%%%%%%%
\n{\bf Example 1.} We want to illustrate Theorem \ref{T-2} on one geometric 
class taken from the list of 93 {\it orbits\,} of the local classification of 
$\Delta^7$ in \cite{dim9}. (The language of \cite{dim9} was different, but the 
objective of that paper was to locally classify $\Delta^7$.) Namely, for 
\,$\C = {\rm GGSTSGS}$ (called 3.1.3.2.3 in \cite{dim9}; it is a single orbit 
in ${\mathbb P}^7\R^2$) we have $s = 2,\,q = 1,\,n_1 = 1$ and $k_0 = l_0 = k_1 = 0$, 
$l_1 = k_2 = 1$, $l_2 = 0$, $l_3 = 2$. By Theorem \ref{T-2}, the table of different 
values appearing in ${\rm der}(\C)$ (consisting of $q + 1 = 2$ rows) is 
%%%%
\begin{align*}
& A_0\,; \\
A_1\Big(&A_0^{+1},\,A_1^{+1},\,A_2^{+1}\Big)\,,
\end{align*}
or else, after the due substitutions, 
%%%%%
\begin{align*}
& 1\,; \\
2\big(&1,\,2,\,5\big)\,.
\end{align*}
By the same Theorem \ref{T-2}, the multiplicities of these values are 
$2 + k_0 + l_0 = 2$, $1 + k_1 + l_1 = 2$, $1 + k_2 + l_2 = 2$, $l_3 - 1 = 1$, 
respectively. Hence ${\rm der}(\C) = (1,\,1,\,2,\,2,\,4,\,4,\,10)$. 
\vskip1mm
\n Let us look now how Jean's recurrences lead to the same result: 
%%%%%%
$$
\begin{array}{rccccccc}
{\rm G} & 1 & & & & & & \\
{\rm G} & 1 & 1 & & & & & \\
{\rm S} & 1 & 1 & 2 & & & & \\
{\rm T} & 1 & 1 & 1 & 3 & & & \\
{\rm S} & 1 & 1 & 2 & 2 & 5 & & \\
{\rm G} & 1 & 1 & 1 & 2 & 2 & 5 & \\
{\rm S} & 1 & 1 & 2 & 2 & 4 & 4 & 10
\end{array}
$$
(omitting the separating commas in the vectors). 
\vskip2.5mm
%%%%%%%%%%%%%%%%%%%%%%%%%R-4
\n{\bf Remark 4.} In \cite{arithm} there is Remark 4.1 which says that, 
in the end of the day, the two-step recurrences from \cite{J} recalled 
in Theorem \ref{com-d}, can be replaced by [much less transparent, however] 
one-step recurrences. We want to underline that it is but an immediate 
consequence of Theorems \ref{T-1} and \ref{T-2}, for the derived vectors 
emerging from these theorems are just {\it functions\,} of the codes of 
geometric classes under consideration. Functions of the extensible set 
of arguments $s,\,k_0,\,l_0,\,k_1,\,l_1,\dots,\,k_s,\,l_s,\,l_{s+1}$. 
\vskip1.2mm
\n To make a point on these one-step recurrences that are now 
virtually at hand: 
%%%%%%%%%%%%%%%%%%%%%%
\begin{itemize}
\item Prolonging a geometric class by G (increasing $l_0$ by 1) 
results in an extremely simple operation on the derived vector. 
\item Prolonging it by T (meaning that $l_0 = 0$) is just increasing 
$k_0$ by 1, with the resulting changes in the derived vector. 
\item Endly, prolonging it by S means: -\,increasing the parameter 
$s$ by 1, -\,the related shift of indices in the parameters $k$ 
and $l$, -\,the consequent changes in the derived vector, all 
of them governed by Theorems \ref{T-1} and/or \ref{T-2}. 
\end{itemize}
%%%%%%%%%%%%%%%%%%%%%
%%%%%%%%%%%%%%%%%%%%%
%%%%%%%%%%%%%%%%%%%%%
\section{Proof of Theorem \ref{T-1}}\label{pfT-1}
The proof is by induction on $s \ge 0$. For $s = 0$, i.\,e., for the geometric 
classes with just one letter S in their codes ($s + 1 = 1$), the derived vectors 
can be computed directly from Theorem \ref{com-d}. Indeed, for the classes 
${\rm G}_{l_1}{\rm S\,T}_{k_0}{\rm G}_{l_0}$, where $l_1 \ge 2$, we are to 
perform $l_1 - 2 \ge 0$ operations G starting from the function $d^{\,2}$, 
then one operation S, then $k_0 \ge 0$ operations T, and eventually $l_0 \ge 0$ 
operations G. 

We know from the beginning that $d^{\,1} = (1)$, $d^{\,2} = (1,\,1)$. 
So, after $l_1 - 2$ operations G, 
$$
d^{\,l_1} = \big(\underbrace{1,\,1,\dots,\,1}_{l_1}\big) = \big(1_{l_1}\big)\,.
$$
Now, by the Fibonacci-like rule in the operation S, 
%%%
\be\label{1:1}
d^{\,l_1 + 1} = {\rm S}(d^{\,l_1 - 1},\,d^{\,l_1}) = \big(1,\,1,\,2_{l_1 - 1}\big)\,.
\ee
When $k_0 > 0$, we next perform $k_0$ operations T, that is, compute $k_0$ 
subsequent terms in arithmetic progressions with the two initial terms either 
(1,\,1) or (1,\,2)\,:
$$
d^{\,l_1 + 2} = {\rm T}(d^{\,l_1},\,d^{\,l_1 + 1}) = \big(1,\,1,\,1,\,3_{l_1 - 1}\big)
$$
$$
\dots\quad\dots\quad\dots\quad\dots
$$
%%%
\be\label{1:2}
d^{\,l_1 + 1 + k_0} = \big(1_{2 + k_0},\,(2 + k_0)_{l_1 - 1}\big)\,.
\ee
Observe that, when $k_0 = 0$, the expression on the right in (\ref{1:1}) 
coincides with the expression on the right in (\ref{1:2}). So (\ref{1:2}) 
is also valid for $k_0 = 0$. Now, applying the operation G $l_0$ times 
to the sequence on the RHS of (\ref{1:2}), we get 
$$
d^{\,l_1 + 1 + k_0 + l_0} = \big(1_{2 + k_0 + l_0},\,(2 + k_0)_{l_1 - 1}\big) = 
\Big(\underbrace{A_0,\dots,\,A_0}_{2 + k_0 + l_0},\,\underbrace{A_1,\dots,\,A_1}_{l_1 - 1}\Big), 
$$
as stated in theorem's statement, irrespectively of whether $l_0 > 0$ or 
$l_0 = 0$. The beginning of induction is done. 
\vskip2mm
%%%%%%%%%%%%%%%%%%%%%%%%%%%%%%%
We are now to justify, for any fixed integer value $s \ge 0$, the induction 
step `$s \Rightarrow s + 1$'. That is, we assume theorem's statement for 
classes with {\bf not more than} $s + 1$ letters S in their codes and work 
towards the analogous statement for an arbitrary class $\C$ having $s + 2$ 
letters S and with such discrete parameters that $l_1 = l_2 = \dots = l_{s+1} = 0$. 

It is natural to truncate the code of $\C$ just before the last from the 
left (the first from the right!) letter S in it. The inductive hypothesis 
will be applied to this truncated word, with the up-shift of indices in 
theorem's parameters by 1 and -- attention -- with the last parameter $l_0$ 
in the very wording of theorem (its role played now by $l_1$) vanishing. 
In parallel, the hypothesis will also be applied to this truncated word 
with {\it its\,} last (or: first from the right) letter {\it deleted}. 
We are going to treat their derived vectors as known and explicitly 
described by the theorem under proof. 

\n Observe that the latter (doubly truncated) word, may have either $s + 1$ 
or $s$ letters S. The first case occurs, clearly, when, there is a gap in 
the code of $\C$ between its $(s + 1)$-st and $(s + 2)$-nd letter S counted 
from the left (between the first and second letter S counted from the right). 
That is, when $k_1 \ge 1$. The second case -- when there is no mentioned gap, 
or else when $k_1 = 0$. We will argue completely separately in each case, 
for the information concerning the derived vector of the doubly truncated 
word depends critically on the case occurring. 
%%%%%%%%%%%%%%%%%%
%%%%%%%%%%%%%%%%%%
\subsection{Case $k_1 \ge 1$.}
In this case, the above-mentioned truncated words are of the form 
${\bf W}{\rm ST}_{k_1}$ and ${\bf W}{\rm ST}_{k_1 - 1}$, where {\bf W} is 
an admissible word with $s$ letters S. 

\n So, by the inductive assumption, the vector ${\rm der}({\bf W}{\rm ST}_{k_1})$ 
is non-decreasing and features (with due multiplicities precised in the 
inductive assumption) the numbers 
$$
A_0^{+1},\ A_1^{+1},\,\dots,\ A_s^{+1},\ A_{s+1}^{+1},
$$
while the vector ${\rm der}({\bf W}{\rm ST}_{k_1 - 1})$ is also non-decreasing 
and features formally `similar' numbers (with their proper multiplicities 
also precised in the inductive assumption), yet built on the basis of 
the initial data $k_1 - 1,\,k_2,\dots,\,k_{s+1}$ instead of 
$k_1,\,k_2,\dots,\,k_{s+1}$ in the previous case. 
Let us write those numbers down as  
$$
\O{A}_0^{+1},\ \O{A}_1^{+1},\,\dots,\ \O{A}_s^{+1},\ \O{A}_{s+1}^{+1}.
$$
These are our basic construction bricks. With their help we will 
firstly express the vector 
%%%%
\be\label{1:3}
{\rm S\Big({\rm der}(\bf W}{\rm ST}_{k_1 - 1}),\,\,{\rm der}\big({\bf W}{\rm ST}_{k_1}\big)\,\Big)\,.
\ee
In order to perform this operation S, one writes its vector 
arguments one above the other, identically indented on the left. 
\vskip1.4mm
\n So, in one row, called ($*$), from the left to right: 
%%%%%%%%%%%
\begin{itemize}
\item $2 + k_1 - 1 = 1 + k_1$ times $\O{A}_0^{+1}$, 
then 
\item $1 + k_2$ times $\O{A}_1^{+1}$, then 
\item $\quad\dots\quad\dots\quad\dots\quad\dots\quad\dots$
\item $1 + k_{s+1}$ times $\O{A}_s^{+1}$, and then 
\item $l_{s+2} - 1$ times $\O{A}_{s+1}^{+1}$.
\end{itemize}

\n And then in the following row, called ($**$), 
also from the left to right: 
%%%%%%%%%
\begin{itemize}
\item $2 + k_1$ times $A_0^{+1}$, then 
\item $1 + k_2$ times $A_1^{+1}$, then 
\item $\quad\dots\quad\dots\quad\dots\quad\dots$
\item $1 + k_{s+1}$ times $A_s^{+1}$, and, eventually, 
\item $l_{s+2} - 1$ times $A_{s+1}^{+1}$.
\end{itemize}

\n In the subsequent row, indented on the left exactly 
as ($*$) and ($**$) are, we ought, by the definition of S, 
to start with two 1's and then perform the Fibonacci rule 
skew-wards from NW to SE. Let us present this procedure 
in detail in the initial parts of the row ($*$) (with 
its $1 + k_1$ initial entries) and ($**$) (with its 
$2 + k_1$ initial entries):
$$
\begin{array}{cccccc}
\O{A}_0^{+1} & \dots    & \O{A}_0^{+1}            & \O{A}_0^{+1} &                         &  \\
A_0^{+1}     & A_0^{+1} & \dots                   & A_0^{+1}     & A_0^{+1}                &  \\
{\bf 1}      & {\bf 1}  & A_0^{+1} + \O{A}_0^{+1} & \dots        & A_0^{+1} + \O{A}_0^{+1} & A_0^{+1} + \O{A}_0^{+1}
\end{array}
$$
Given that we know the entire rows ($*$) and ($**$), it 
is now visible that the overall result of the operation S 
is a row, called ($***$), by one entry longer than ($**$), 
in which there go, from the left to right, 
%%%%%%%%%%%
\begin{itemize}
\item 2 times 1, then 
\item $1 + k_1$ times $A_0^{+1} + \O{A}_0^{+1}$, then 
\item $1 + k_2$ times $A_1^{+1} + \O{A}_1^{+1}$, then 
\item $\quad\dots\quad\dots\quad\dots\quad\dots\quad\dots$
\item $1 + k_{s+1}$ times $A_s^{+1} + \O{A}_s^{+1}$, and then 
\item $l_{s+2} - 1$ times $A_{s+1}^{+1} + \O{A}_{s+1}^{+1}$.
\end{itemize}

Now we are to perform the operation T certain ($k_0 \ge 0$) 
number of times. When $k_0 > 0$, we apply it for the first 
time to the rows ($**$) and ($***$) (these two are its vector 
arguments). Let us present this procedure in detail, like 
with the preceding operation, in the initial parts of these 
rows, focusing on the $2 + k_1$ initial entries in the row 
($**$) and on the $3 + k_1$ initial entries in the row ($***$):
$$
\begin{array}{ccccccc}
1 & A_0^{+1} & A_0^{+1} & \dots & A_0^{+1} &  &  \\
1 & 1 & A_0^{+1} + \O{A}_0^{+1} & A_0^{+1} + \O{A}_0^{+1} & \dots & A_0^{+1} + \O{A}_0^{+1} & \\
{\bf 1} & {\bf 1} & 1 & A_0^{+1} + 2\O{A}_0^{+1} & A_0^{+1} + 2\O{A}_0^{+1} & \dots & A_0^{+1} + 2\O{A}_0^{+1}
\end{array}
$$
Given that we know the entire rows ($**$) and ($***$), it 
now becomes clear (a)\,how the rest of the output of this 
first operation T looks like, and (b)\,how the entire 
continuation with the consecutive $k_0 - 1$ operations 
T looks like. Namely, the resulting eventual derived 
vector consists of the entries 
%%%%%%%%%
\begin{itemize}
\item $2 + k_0$ times $1 \,(=\,A_0)$, then 
\item $1 + k_1$ times $A_0^{+1} + (1 + k_0)\O{A}_0^{+1}$, then 
\item $1 + k_2$ times $A_1^{+1} + (1 + k_0)\O{A}_1^{+1}$, then
\item $\quad\dots\quad\dots\quad\dots\quad\dots\quad\dots$
\item $1 + k_{s+1}$ times $A_s^{+1} + (1 + k_0)\O{A}_s^{+1}$, 
and then
\item $l_{s+2} - 1$ times $A_{s+1}^{+1} + (1 + k_0)\O{A}_{s+1}^{+1}$.
\end{itemize}
Observe that, when $k_0 = 0$ (i.\,e., when there is no operation T 
after the S) this result also takes effect, because it coincides 
with the row ($***$). 
\vskip1.5mm
What we still need in this case $k_1 \ge 1$ is 
%%%%%%%%%
\begin{lem}\label{key}
$A_j^{+1} + (1 + k_0)\O{A}_j^{+1} \,= \,A_{j+1}$ \,\,\,for \,\,\,$j = 0,\,1,\dots,\,s,\,s+1$. 
\end{lem}
Proof by induction on $j$. The beginning of induction, for $j = 0$ 
and $j = 1$, is immediate. Indeed, $A_0^{+1} + (1 + k_0)\O{A}_0^{+1} 
= 1 + (1 + k_0)\cdot 1 = A_1$ and 
$$
A_1^{+1} + (1 + k_0)\O{A}_1^{+1} = 2 + k_1 + (1 + k_0)(2 + k_1 - 1) 
= 1 + (2 + k_0)(1 + k_1) = A_0 + (1 + k_1)A_1 = A_2\,.
$$
The induction step, in the form $(j-1,\,j) \Rightarrow j + 1$ 
for any fixed $j$, $1 \le j \le s$, is only a bit longer. 
From the definition of the $A$ sequences 
%%%
\be\label{def_1}
A_{j+1}^{+1} = A_{j-1}^{+1} + (1 + k_{j+1})A_j^{+1},
\ee
%%%%%%
\be\label{def_2}
\O{A}_{j+1}^{+1} = \O{A}_{j-1}^{+1} + (1 + k_{j+1})\O{A}_j^{+1}.
\ee
Now we check if lemma's statement holds for $j + 1$: 
%%%%%%
\begin{align*}
A_{j+1}^{+1} &+ (1 + k_0)\O{A}_{j+1}^{+1} \stackrel{\text{by (\ref{def_1}) and (\ref{def_2})}}{=} 
A_{j-1}^{+1} + (1 + k_{j+1})A_j^{+1} + (1 + k_0)\big(\O{A}_{j-1}^{+1} + (1 + k_{j+1})\O{A}_j^{+1}\big)\\
&= A_{j-1}^{+1} + (1 + k_0)\O{A}_{j-1}^{+1} + (1 + k_{j+1})\big(A_j^{+1} + (1 + k_0)\O{A}_j^{+1}\big)\\
&\ \stackrel{\text{by ind. ass.}}{=} A_j + (1 + k_{j+1})A_{j+1} = A_{j+2}\,.
\end{align*}
Taking into account this lemma, after the operations S and $k_0 \ge 0$ times T 
we have the result row in a clearer form 
%%%%%%%%%
\begin{itemize}
\item $2 + k_0$ times $A_0$, then 
\item $1 + k_1$ times $A_1$, then 
\item $1 + k_2$ times $A_2$, then
\item $\quad\dots\quad\dots\quad\dots\quad\dots$
\item $1 + k_{s+1}$ times $A_{s+1}$, 
and then
\item $l_{s+2} - 1$ times $A_{s+2}$.
\end{itemize}
This is the vector ${\rm der}\big({\bf W}{\rm ST}_{k_1}{\rm ST}_{k_0}\big)$. 
To finish the case $k_1 \ge 1$, we remember that there still are 
$l_0 \ge 0$ letters G in the end of the code of $\C$. But the operation 
G is extremely simple: one just adds, from the left, one entry 1. 
Hence the eventual vector 
${\rm der}\big({\bf W}{\rm ST}_{k_1}{\rm ST}_{k_0}{\rm G}_{l_0}\big)$ 
has just $l_0$ entries 1 ($=\,A_0$) more in the beginning. 
That is, it consists of 
%%%%%%%%%
\begin{itemize}
\item $2 + k_0 + l_0$ times $A_0$, then 
\item $1 + k_1$ times $A_1$, then 
\item $1 + k_2$ times $A_2$, then
\item $\quad\dots\quad\dots\quad\dots\quad\dots$
\item $1 + k_{s+1}$ times $A_{s+1}$, 
and then
\item $l_{s+2} - 1$ times $A_{s+2}$.
\end{itemize}
Case $k_1 \ge 1$ [in the inductive proof of Theorem \ref{T-1}] is now proved. 
%%%%%%%%%%%%%%%%%%%%%
%%%%%%%%%%%%%%%%%%%%%
\subsection{Case $k_1 = 0$.}

Now the two truncated words are just of the form {\bf W}S and {\bf W}, 
where {\bf W} is an admissible word with $s$ letters S. 

\n The word {\bf W}S is admissible, too, and has $s + 1$ letters S. 
So, by the inductive assumption, the vector der({\bf W}S) is non-decreasing 
and features, with due multiplicities precised in the inductive 
assumption, the numbers, written in the growing order: 
$$
A_0^{+1},\ A_1^{+1},\,\dots,\ A_s^{+1},\ A_{s+1}^{+1}\,,
$$
now with an important restriction $k_1 = 0$, hence $A_1^{+1} = 2 + k_1 = 2$ 
etc. The word {\bf W} having $s$ letters S, by the same assumption, 
the vector der({\bf W}) is non-decreasing and features the numbers, 
in the growing order (with their proper multiplicities 
also precised in the inductive assumption), 
$$
A_0^{+2},\ A_1^{+2},\,\dots,\ A_s^{+2}\,.
$$
These are our construction bricks this time. We first have to compute 
the derived vector S\big(der({\bf W}),\,der({\bf W}S)\big). To this 
end we expand der({\bf W}) in a row, called ($*'$), consisting of 
%%%%%%%%%
\begin{itemize}
\item $2 + k_2$ times $A_0^{+2}$, then 
\item $1 + k_3$ times $A_1^{+2}$, then 
\item $\quad\dots\quad\dots\quad\dots\quad\dots$
\item $1 + k_{s+1}$ times $A_{s-1}^{+2}$, and, eventually, 
\item $l_{s+2} - 1$ times $A_s^{+2}$. 
\end{itemize}
(Clearly, when $s = 0$, then only the last group is present here 
and der({\bf W}) consists of $l_2 - 1$ values $A_0^{+2} = 1$.)  

\n And we also expand der({\bf W}S), in a row called ($**'$), 
only slightly differing from the row ($**$) in the previous 
case. The difference resides in the current specification 
$k_1 = 0$ (which generates, naturally, a different sequence 
of values $A$ than when $k_1 > 0$)\,:
%%%%%%%%%
\begin{itemize}
\item 2 times $A_0^{+1}$, then 
\item $1 + k_2$ times $A_1^{+1}$, then 
\item $\quad\dots\quad\dots\quad\dots\quad\dots$
\item $1 + k_{s+1}$ times $A_s^{+1}$, and, eventually, 
\item $l_{s+2} - 1$ times $A_{s+1}^{+1}$.
\end{itemize}
The continuation is, basically, already known. The row ($*'$) 
is written above row ($**'$), both identically indented on the 
left, and in the following row below them, called ($***'$) and 
identically indented on the left, we start with two 1's written 
in bold to put them in relief, then apply the Fibonacci-like rule 
from the NW to SE. Let us trace down this procedure in the initial 
parts of the rows: first $2 + k_2$ terms in ($*'$) and first 
$3 + k_2$ terms in ($**'$): 
$$
\begin{array}{rcccccc}
(*')   & 1 & A_0^{+2} & \dots & A_0^{+2} &  &  \\
(**')  & 1 & 1 & A_1^{+1} & \dots & A_1^{+1} &  \\
(***') & {\bf 1} & {\bf 1} & 2 & A_1^{+1} + A_0^{+2} & \dots & A_1^{+1} + A_0^{+2} 
\end{array}
$$
(the first entry in ($*'$) is written as 1, and not $A_0^{+2}$). 
Given the above-listed multiplicities of different values in the 
rows ($*'$) and ($**'$), the pattern observed on the initial terms 
of both the arguments and result smoothly extends to the entire 
rows. The outcome row $(***') = {\rm S}\big((*'),\,(**')\big)$ 
consists of entries, in the order of their appearing: 
%%%%%%%%%%%
\begin{itemize}
\item 2 times 1, then 
\item 1 time 2, then 
\item $1 + k_2$ times $A_1^{+1} + A_0^{+2}$, then 
\item $1 + k_3$ times $A_2^{+1} + A_1^{+2}$, then 
\item $\quad\dots\quad\dots\quad\dots\quad\dots\quad\dots$
\item $1 + k_{s+1}$ times $A_s^{+1} + A_{s-1}^{+2}$, and then 
\item $l_{s+2} - 1$ times $A_{s+1}^{+1} + A_s^{+2}$.
\end{itemize}
Now we are to continue with the operation T performed $k_0$ times. 
Since T is a (skew-wards NW to SE) arithmetic-progression-like 
operation which is initially being applied after the operation S 
acting on the rows: ($*'$) as the first argument in S and ($**'$) 
as second argument, the arithmetic progressions obtained in the 
outcome of $k_0$ operations T are either constant and equal to 1, 
or have their {\it differences\,} in the row ($*'$) -- just entries 
in that row. In fact, these are the sequences 
%%%%%%%%%
\begin{itemize}
\item $1,\,1,\,1,\dots,\,1$ (altogether $k_0$ sequences of lengths $3,\,4,\dots,\,2 + k_0$), 
\item $1,\,2,\,3,\dots,\,2 + k_0$ (1 time), 
\item $A_1^{+1},\,\,A_1^{+1} + A_0^{+2},\,\,A_1^{+1} + 2A_0^{+2},\dots,\,\,A_1^{+1} + (1 + k_0)A_0^{+2}$ ($1 + k_2$ times), 
\item $A_2^{+1},\,\,A_2^{+1} + A_1^{+2},\,\,A_2^{+1} + 2A_1^{+2},\dots,\,\,A_2^{+1} + (1 + k_0)A_1^{+2}$ ($1 + k_3$ times),
\item $\qquad\dots\qquad\dots\qquad\dots\qquad\dots\qquad\dots$
\item $A_s^{+1},\,\,A_s^{+1} + A_{s-1}^{+2},\,\,A_1^{+1} + 2A_{s-1}^{+2},\dots,\,\,A_1^{+1} + (1 + k_0)A_{s-1}^{+2}$ ($1 + k_{s+1}$ times),
\item $A_{s+1}^{+1},\,\,A_{s+1}^{+1} + A_s^{+2},\,\,A_{s+1}^{+1} + 2A_s^{+2},\dots,\,\,A_{s+1}^{+1} + (1 + k_0)A_s^{+2}$ ($l_{s+2} - 1$ times)
\end{itemize}
(except for the irregular sequences in the beginning, the arithmetic 
sequences here have length $2 + k_0$, and their first (resp., second) 
terms are always in the row ($**'$) (resp., ($***'$)\,). 
\vskip1mm
\n Summing up, the vector ${\rm der}\big({\bf W}{\rm SST}_{k_0}\big)$ 
is non-decreasing and consists of 
%%%%%%%%%%%%%
\begin{itemize}
\item $2 + k_0$ times $1 =\,A_0$, then
\item 1 time $2 + k_0 = A_1$, then 
\item $1 + k_2$ times $A_1^{+1} + (1 + k_0)A_0^{+2}$, then 
\item $1 + k_3$ times $A_2^{+1} + (1 + k_0)A_1^{+2}$, then 
\item $\quad\dots\quad\dots\quad\dots\quad\dots\quad\dots$ 
\item $1 + k_{s+1}$ times $A_s^{+1} + (1 + k_0)A_{s-1}^{+2}$, 
and eventually 
\item $l_{s+2} - 1$ times $A_{s+1}^{+1} + (1 + k_0)A_s^{+2}$.
\end{itemize}
To proceed, analogously as in the previous case $k_1 > 0$, 
we need 
%%%%%%%%%%%
\begin{lem}\label{key_bis}
When \,$k_1 = 0$, then \,$A_j^{+1} + (1 + k_0)A_{j-1}^{+2} = A_{j+1}$ 

\n for \,$j = 1,\,2,\dots,\,s,\,s+1$.
\end{lem}
{\it Attention.} There is no use to directly compare this lemma with 
Lemma \ref{key}, because now we are under the restriction $k_1 = 0$. 
The two lemmas deal with {\it disjoint\,} families of integer 
sequences $A$. 
\vskip1.5mm
%%%%%%%%%%%%%%%%%%%%%%%%
\n Proof of lemma is inductive on $j$. For $j = 1$ it is quick, 
$A_1^{+1} + (1 + k_0)A_0^{+2} = 2 + k_1 + 1 + k_0 = 1 + 2 + k_0 
= A_0 + A_1 = A_0 + (1 + k_1)A_1 = A_2$. 

\n For $j = 2$ it is a bit longer, 
%%%%%%%%%%
\begin{align*}
A_2^{+1} + (1 + k_0)A_1^{+2} &= 1 + (1 + k_2)A_1^{+1} + (1 + k_0)(2 + k_2) \\
&= 1 + 2(1 + k_2) + (1 + k_0)(1 + 1 + k_2) = 2 + k_0 + (1 + k_2)(1 + 2 + k_0)\\
&= A_1 + (1 + k_2)(A_0 + A_1) = A_1 + (1 + k_2)(A_0 + (1 + k_1)A_1) \\
&= A_1 + (1 + k_2)A_2 = A_3\,.
\end{align*}
As for the inductive step, again in the form $(j-1,\,j) \Rightarrow j + 1$, 
with, this time, $j \in \{2,\,3,\dots,\,s\}$, we will nearly mimick the 
proof of Lemma \ref{key}. We will again make use of identity (\ref{def_1}) 
and of a similar identity 
%%%
\be\label{def_3}
A_j^{+2} = A_{j-2}^{+2} + (1 + k_{j+1})A_{j-1}^{+2}
\ee
which follows from the definition of the $A$ sequences as well. 
\vskip1.5mm
\n Let us check now if lemma's statement holds for $j + 1$: 
%%%%%%
\begin{align*}
A_{j+1}^{+1} &+ (1 + k_0)A_j^{+2} \stackrel{\text{by (\ref{def_1}) and (\ref{def_3})}}{=} 
A_{j-1}^{+1} + (1 + k_{j+1})A_j^{+1} + (1 + k_0)\big(A_{j-2}^{+2} + (1 + k_{j+1})A_{j-1}^{+2}\big)\\
&= A_{j-1}^{+1} + (1 + k_0)A_{j-2}^{+2} + (1 + k_{j+1})\big(A_j^{+1} + (1 + k_0)A_{j-1}^{+2}\big)\\
&\ \stackrel{\text{by ind. ass.}}{=} A_j + (1 + k_{j+1})A_{j+1} = A_{j+2}\,.
\end{align*}
Lemma \ref{key_bis} is now proved by induction. 
\vskip2.5mm
We are about to finish the proof of the case $k_1 = 0$ in the theorem. 
Indeed, Lemma \ref{key_bis} simplifies the expressions for the entries 
in the vector ${\rm der}\big({\bf W\rm SST}_{k_0}\big)$, while the $l_0$ 
operations G that are still to be performed add $l_0$ entries 1 on the 
left. Hence the vector ${\rm der}\big({\bf W\rm SST}_{k_0}{\rm G}_{l_0}\big)$, 
is non-decreasing and consists of 
%%%%%%%%%%%%%
\begin{itemize}
\item $2 + k_0 + l_0$ times $A_0$, then
\item $1 + k_1 = 1$ time $A_1$, then 
\item $1 + k_2$ times $A_2$, then 
\item $1 + k_3$ times $A_3$, then 
\item $\quad\dots\quad\dots\quad\dots$ 
\item $1 + k_{s+1}$ times $A_{s+1}$, 
and eventually 
\item $l_{s+2} - 1$ times $A_{s+2}$.
\end{itemize}
Case $k_1 = 0$ [in the inductive proof of Theorem \ref{T-1}] is finished. 
\vskip1.5mm
\n Theorem \ref{T-1} is proved. 
%%%%%%%%%%%%%%%%%%
%%%%%%%%%%%%%%%%%%
\subsection{Reformulations of Theorem \ref{T-1}.}
In order to simplify our to-be-done work in proving Theorem \ref{T-2} in 
Section \ref{pfT-2}, let us slightly reformulate Theorem \ref{T-1}, then 
draw a corollary from that reformulation, and then yet another, more 
far-reaching, corollary. 
%%%%%%%%%%%%
\begin{obs}\label{reform}
In the notation from {\rm \,Theorem \ref{T-1}}, for an admissible word with 
$s + 1$ letters {\rm \,S} and such that \,$l_1 = l_2 = \cdots = l_s = 0$, 
the sequence of operations 
$$
{\rm S},\:{\rm T}^{k_s},\:{\rm S},\:{\rm T}^{k_{s-1}},\,\dots,\,{\rm S},\:{\rm T}^{k_1},\:
{\rm S},\:{\rm T}^{k_0},\:{\rm G}^{l_0}
$$
started on the pair of arguments $\big(\underbrace{1,\,1,\dots,\,1}_{l_{s+1} - 1}\big)$ 
and $\big(\underbrace{1,\,1,\dots,\,1}_{l_{s+1}}\big)$ produces as the eventual result 
the vector described in the statement of {\rm \,Theorem \ref{T-1}}. 
\end{obs}
This is just a reformulation of Theorem \ref{T-1}, compare the beginning of its proof. 
%%%%%%%%%%%%%
\begin{cor}\label{reformm}
The sequence of operations 
$$
{\rm S},\:{\rm T}^{k_{n-1}},\:{\rm S},\:{\rm T}^{k_{n-2}},\,\dots,\,{\rm S},\:{\rm T}^{k_1},\:
{\rm S},\:{\rm T}^{k_0},\:{\rm G}^{l_0}
$$
($n$ is a natural integer not exceeding $s$) started on the pair of vector arguments 
$\big(\underbrace{w_0,\,w_0,\dots,\,w_0}_{m_0}\big)$ \,and \,$\big(1,\,\underbrace{w_0,\,w_0,\dots,\,w_0}_{m_0}\big)$, 
where $w_0$ is a positive integer, produces in the result the vector with the entries, 
consecutively: 
\end{cor}
%%%%%%%%%
\begin{itemize}
\item $2 + k_0 + l_0$ times $A_0$, then 
\item $1 + k_1$ times $A_1$, then 
\item $\quad\dots\quad\dots\quad\dots$
\item $1 + k_{n-1}$ times $A_{n-1}$, and eventually 
\item $m_0$ times $A_n w_0$\,.
\end{itemize}
This corollary is only slightly more general than Observation \ref{reform} 
(instead of $s$ there is $n - 1$ now, $m_0$ instead of $l_{s+1} - 1$, 
and $w_0$ instead of 1) and its method of proof is the same as that 
in Theorem \ref{T-1}. 
%%%%%%%%%%%
\begin{cor}\label{reformmm}
The same sequence of operations as in {\rm \,Corollary \ref{reformm}}, when 
started on the pair of longer vector arguments 
$$
\big(\underbrace{w_0,\,w_0,\dots,\,w_0}_{m_0},\:\underbrace{w_1,\,w_1,\dots,\,w_1}_{m_1},\:\dots,\:\underbrace{w_N,\,w_N,\dots,\,w_N}_{m_N}\big)
$$
and 
$$
\big(1,\:\underbrace{w_0,\,w_0,\dots,\,w_0}_{m_0},\:\underbrace{w_1,\,w_1,\dots,\,w_1}_{m_1},\:\dots,\:\underbrace{w_N,\,w_N,\dots,\,w_N}_{m_N}\big)\,,
$$
where $w_0,\,w_1,\dots,\,w_N$ are positive integers, produces in the result 
the vector with the entries, consecutively: 
\end{cor}
%%%%%%%%%
\begin{itemize}
\item $2 + k_0 + l_0$ times $A_0$, then 
\item $1 + k_1$ times $A_1$, then 
\item $\quad\dots\quad\dots\quad\dots\quad\dots$
\item $1 + k_{n-1}$ times $A_{n-1}$, then 
\item $m_0$ times $A_n w_0$, then 
\item $m_1$ times $A_n w_1$, then 
\item $\quad\dots\quad\dots\quad\dots$
\item $m_N$ times $A_n w_N$\,.
\end{itemize}
This is a direct consequence of Corollary \ref{reformm}. Indeed, replacing 
in that corollary $w_0$ and $m_0$ by $w_i$ and $m_i$, respectively, $i \in \{1,\,2,\dots,\,N\}$, 
and then {\it concatenating\,} the vector results for $i = 1,\,2,\dots,\,N$ 
with the initial vector result in Corollary \ref{reformm}, we get the statement 
in Corollary \ref{reformmm}. 
%%%%%%%%%%%%%%%%%%
%%%%%%%%%%%%%%%%%%
%%%%%%%%%%%%%%%%%%
\section{Proof of Theorem \ref{T-2}}\label{pfT-2}
The proof of Theorem \ref{T-2} will be by induction on $q \in \{0,\,1,\dots,\,s\}$ 
in all the pairs of integers $(q,\,s)$, $0 \le q \le s$. The size of $s \ge q$ 
will not play any role in the arguments.
\vskip1.5mm
\n Before doing this we want to formulate and prove a statement that would greatly 
simplify the inductive step in that announced induction. Note that the notation used 
in Proposition \ref{P-2} below is fully compatible with that of Corollary \ref{reformmm}. 
It is done purposefully -- in order for the reader to see that it is basically 
Corollary \ref{reformmm} again, only given in other terms.
%%%%%%%%%%%%%%%%%%%%%%%%
\begin{prop}\label{P-2}
Let the code of a geometric class \,$\C$ have the form 
$$
{\bf W}{\rm GST}_{k_{n-1}}{\rm ST}_{k_{n-2}}\dots{\rm ST}_{k_1}{\rm ST}_{k_0}{\rm G}_{l_0}\,,
$$
{\bf W} -- the code of an ancestor of \,$\C$, ending on {\rm \,S} or {\rm \,T}, 
whose derived vector {\rm der}$({\bf W})$ is non-decreasing and has the different 
values $1 = w_0 < w_1 < \cdots < w_N$ appearing in it with multiplicities 
$m_0,\,m_1,\dots,\,m_N$, respectively. Then ${\rm der}(\C)$ is non-decreasing 
as well and consists of 
\end{prop}
%%%%%
\begin{itemize}
\item $2 + k_0 + l_0$ times $A_0$, then 
\item $1 + k_1$ times $A_1$, then 
\item $\quad\dots\quad\dots\quad\dots\quad\dots$ 
\item $1 + k_{n-1}$ times $A_{n-1}$, then 
\item $m_0$ times $A_n\,(\,= A_n w_0)$, then 
\item $m_1$ times $A_n w_1$, then 
\item $\quad\dots\quad\dots\quad\dots$
\item $m_N$ times $A_n w_N$. 
\end{itemize}
Proof of Proposition \ref{P-2}. It is a direct consequence of Corollary \ref{reformmm}, 
because the derived vector of {\bf W} can be taken as the first vector argument in 
that corollary, upon which the derived vector of {\bf W}G becomes the second vector 
argument in that corollary. (The operation G adds a single 1 on the left, while the 
particular value $w_0 = 1$ in the proposition does not interfere with this action.) 
\vskip2.5mm
%%%%%%%%%%%%%%%%%%%%%%%%%%%
{\bf The beginning of induction} for $q = 0$ {\bf in Theorem \ref{T-2}} is already 
done -- it is Theorem \ref{T-1}, proved already for all $s \ge 0$. 
\vskip2mm
%%%%%%%%%%%%%%
{\bf The inductive step} `$q - 1 \Rightarrow q$' {\bf in Theorem \ref{T-2}} for any 
fixed $q \ge 1$ and arbitrary $s \ge q$, assuming the due information for $q - 1$ 
and all $s' \ge q - 1$. 
\vskip1mm
\n Let us focus our attention on the one before last from the left (or last, if $l_0 = 0$) 
segment of G's in the word of \,$\C$ and truncate this word, calling {\bf W} the result 
of truncation, {\it before\,} the last letter G in that segment. That is, take 
$$
{\bf W} \,= \,{\rm G}_{l_{s+1}}\!{\rm S\,T}_{k_s}\dots{\rm S\,T}_{k_{n_1}}\!{\rm G}_{l_{n_1} - 1}\, 
$$
where $l_{n_1} - 1 \ge 0$ (the value 0 is not excluded). Then the vector der({\bf W}) 
is known from the inductive assumption, for there are only $q - 1$ positive values 
among the parameters $l_{n_1 + 1},\,l_{n_1 + 2},\dots,\,l_s$ related to the word 
{\bf W}, and $s' = s - n_1 \ge q - 1$. (Indeed, from our construction, $s - n_1 \ge 
n_q - n_1 > n_{q-1} - n_1 > \cdots > n_2 - n_1 \ge 1$.) So what is precisely der({\bf W})? 
\vskip1.5mm
The inductive assumption {\it formally\,} would have yielded a description in terms 
of the set of letters (as if a formal alphabet) 
$$
\{n_1,\,n_2,\dots,\,n_{q-1},\,s - n_1,\,k_0,\,l_0,\dots,\,k_{s-n_1},\,l_{s-n_1},\,l_{s - n_1 + 1}\}\,, 
$$
because of the relevant bare cardinalities related to the word {\bf W}. Yet, in the concrete 
case of {\bf W}, all {\it indices\,} in the parameters $k$ and $l$ are naturally raised by 
$+\,n_1$, ranging from $n_1$ to $s + 1$ (and not from 0 to $s - n_1 + 1$). And the {\it distances\,} 
in indices are $n_2 - n_1,\,\,n_3 - n_1,\,\,\dots,\,\,n_q - n_1$ instead of $n_1 - 0,\,\,n_2 - 0,\,
\dots,\,\,n_{q-1} - 0$, respectively (the distances in indices in {\bf W} are from $n_1$, 
not from 0). 
\vskip1mm
\n Taking into account all these specific features, the data about {\bf W} which follow 
from the inductive assumption, and are about to be plugged into Proposition \ref{P-2}, 
read: 
\vskip1mm
In the vector der({\bf W}) there appear $s - n_1 + 2$ different values 
$1 = w_0 < w_1 < w_2 < \dots < w_{s - n_1 + 1}$, listed below in the 
growing order, in the following $(q - 1) + 1 = q$ separate rows: 
%%%%%%%%%%%
\begin{align*}
& A_0^{+n_1},\;A_1^{+n_1},\,\dots,\;A_{n_2 - n_1 - 1}^{+n_1}\,; \\
%%%
A_{n_2 - n_1}^{+n_1}\big( & A_0^{+n_2},\;A_1^{+n_2},\,\dots,\;A_{n_3 - n_2 - 1}^{+n_2}\big)\,; \\
%%%
\dots \qquad & \dots \qquad \dots \qquad \dots \qquad \dots\\
%%%
\prod_{j=1}^{q-2}A_{n_{j+1} - n_j}^{+n_j}\Big( & A_0^{+n_{q-1}},\;A_1^{+n_{q-1}},\,\dots,
\;A_{n_q - n_{q-1} - 1}^{+n_{q-1}}\Big)\,; \\
%%%
\prod_{j=1}^{q-1}A_{n_{j+1} - n_j}^{+n_j}\Big( & A_0^{+n_q},\;A_1^{+n_q},\,\dots,
\;A_{s - n_q + 1}^{+n_q}\Big)\,.
\end{align*}
Moreover, der({\bf W}) is non-decreasing and the multiplicities of its above-listed 
different values are as follows. 
%%%%%%%
\begin{itemize}
\item $w_0 = A_0^{+n_1}$ appears $m_0 = 2 + k_{n_1} + (l_{n_1} - 1) = 1 + k_{n_1} + l_{n_1}$ 
times in row (remember that the last group of G's in {\bf W} consists, by construction, of 
$l_{n_1} - 1$ letters); 
\item for $1 \le j \le s - n_1$, $w_j$ appears $m_j = 1 + k_{n_1 + j} + l_{n_1 + j}$ times in row 
(that is, $1 + k_{n_1 + j}$ times, when $j \notin \{n_2 - n_1,\,n_3 - n_1,\,\dots,\,n_q - n_1\}$);
\item the last and biggest $w_{s - n_1 + 1}$ appears $m_{s - n_1 + 1} = l_{s+1} - 1$ times in row. 
\end{itemize}
\vskip1mm
%%%%%%%%%%%%%%%%%%%
We are about to apply Proposition \ref{P-2} to the code of the class \,$\C$ under the induction 
step procedure and to its truncation {\bf W}, for $n = n_1$, $N = s - n_1 + 1$, the different 
values $w_i$ in der({\bf W}) and their multiplicities $m_i$ ($i = 0,\,1,\dots,\,s - n_1 + 1$), 
as specified above. 
\vskip1.5mm
\n So Proposition \ref{P-2} says that the table of different values in ${\rm der}(\C)$ 
is being obtained by just multiplying by $A_{n_1}$ all the entries in the table of $q$ 
rows for der({\bf W}) above {\it and\,} by inserting on top of them the single row 
$A_0,\,A_1,\dots,\,A_{n_1 - 1}$. That is, the table of different values showing up 
in the vector ${\rm der}(\C)$ is of the form 
%%%%%%%%%%%%
\begin{align*}
&A_0,\,A_1,\,\dots,\,A_{n_1 - 1}\,; \\
%%%
A_{n_1}\big( & A_0^{+n_1},\;A_1^{+n_1},\,\dots,\;A_{n_2 - n_1 - 1}^{+n_1}\big)\,; \\
%%%
A_{n_1}A_{n_2 - n_1}^{+n_1}\big( & A_0^{+n_2},\;A_1^{+n_2},\,\dots,\;A_{n_3 - n_2 - 1}^{+n_2}\big)\,; \\
%%%
\dots \qquad & \dots \qquad \dots \qquad \dots \qquad \dots\\
%%%
A_{n_1}\prod_{j=1}^{q-2}A_{n_{j+1} - n_j}^{+n_j}\Big( & A_0^{+n_{q-1}},\;A_1^{+n_{q-1}},\,\dots,
\;A_{n_q - n_{q-1} - 1}^{+n_{q-1}}\Big)\,; \\
%%%
A_{n_1}\prod_{j=1}^{q-1}A_{n_{j+1} - n_j}^{+n_j}\Big( & A_0^{+n_q},\;A_1^{+n_q},\,\dots,
\;A_{s - n_q + 1}^{+n_q}\Big)\,, 
\end{align*}
precisely as expected. Moreover, also by Propositiopn \ref{P-2}, the vector ${\rm der}(\C)$ 
is non-decreasing and the multiplicities of its listed different values which emerge from 
Proposition \ref{P-2} are identical (recall that $n = n_1$) with those in the statement 
of Theorem \ref{T-2}. Theorem's statements for $q$ positive parameters among 
$\{l_1,\,l_2,\dots,\,l_s\}$ are now justified and the step of induction is 
completed. 
\vskip1.5mm
\n Theorem \ref{T-2} is now proved. 
\vskip3mm
%%%%%%%%%%%%%%%%%%%%%%%%%%%%%%%
\n{\bf Afterword.} In the contribution \cite{MZ2}, among many an issue raised, there 
is a far-reaching Question 9.19, consisting of parts 1, 2, and 3. 
\vskip1.3mm
\n Part 1, especially when understood {\it sensu largo} (is it possible to produce 
Jean's recurrences in pure [Goursat] Monster terms?), is a true challenge. 
\vskip1.3mm
\n Part 2 is elementary and covered by the present work (which, as a matter of record, 
completes \cite{arithm}). 
\vskip1.3mm
\n Part 3 is an invitation to research, for the mapping Pc ({\it Puiseux characteristic\,} 
being associated to a geometric class = a word over $\{$R,\,V,\,T$\}$) is explicit in 
\cite{MZ2}, while the mapping GW ([small] growth vector being associated to a geometric 
class = a word over $\{$G,\,S,\,T$\}$) is recursive in \cite{J}, and explicit in 
\cite{arithm} and in the present contribution. 
%%%%%%%%%%%%%%%%%%
%%%%%%%%%%%%%%%%%%
%%%%%%%%%%%%%%%%%%

\end{document}